\documentclass[notitlepage]{revtex4-1}
\usepackage{graphicx,indentfirst}
\usepackage{amscd}
\usepackage[all]{xy}
\usepackage{amsmath}
\usepackage{amsfonts,amssymb}
\usepackage{amsthm}
\usepackage{latexsym}
\usepackage[american]{babel}
\usepackage{times}
\usepackage{tikz}
\usepackage{verbatim}
\usepackage{hyperref}
\usepackage{etoolbox}
\usetikzlibrary{matrix,arrows,decorations.pathmorphing}
\usepackage[a4paper,left=2.54cm,right=2.54cm,bottom=2.54cm,top=2.54cm]{geometry}
\RequirePackage[T1]{fontenc}

\begin{document}

\title[The Poisson bracket on the family algebra]{The noncommutative Poisson bracket and the deformation of the family algebras}
\author{Zhaoting  Wei}
\affiliation{Department of Mathematics, Indiana University, 831 E 3rd Street, Bloomington, IN 47405, USA}
\email{zhaotwei@indiana.edu}
\date{submitted on April 4, 2015, revised on \today}

\newcommand{\End}{\text{End}}
\newcommand{\ad}{\text{ad}}
\newcommand{\Ad}{\text{Ad}}
\newcommand{\tr}{\text{tr}}
\newcommand{\Pol}{\text{Pol}}
\newcommand{\IPBW}{\text{I}_{\text{PBW}}}
\newcommand{\FPBW}{\text{F}_{\text{PBW}}}
\newcommand{\Hom}{\text{Hom}}
\newcommand{\dH}{\text{d}_{\text{H}}}
\newcommand{\HH}{\text{HH}}
\newcommand{\sll}{\text{sl}}

\newtheorem{thm}{Theorem}[section]
\renewcommand\thethm{\arabic{section}.\arabic{thm}}
\newtheorem{lemma}[thm]{Lemma}
\newtheorem{prop}[thm]{Proposition}
\newtheorem{coro}[thm]{Corollary}
\theoremstyle{definition}\newtheorem{defi}{Definition}
\renewcommand\thedefi{\arabic{section}.\arabic{defi}}
\theoremstyle{remark}\newtheorem{eg}{Example}
\renewcommand\theeg{\arabic{section}.\arabic{eg}}
\theoremstyle{remark}\newtheorem{rmk}{Remark}
\renewcommand\thermk{\arabic{section}.\arabic{rmk}}
\theoremstyle{remark}\newtheorem{ctn}{Caution}
\renewcommand\thectn{\arabic{section}.\arabic{ctn}}

\begin{abstract}
The family algebras are introduced by A.A. Kirillov in 2000. In this paper we study the noncommutative Poisson bracket $P$ on the classical family algebra $\mathcal{C}_{\tau}(\mathfrak{g})$. We show that $P$ controls the first-order $1$-parameter formal deformation from $\mathcal{C}_{\tau}(\mathfrak{g})$ to $\mathcal{Q}_{\tau}(\mathfrak{g})$ where the later is the quantum family algebra. Moreover we will prove that the noncommutative Poisson bracket is in fact a Hochschild $2$-coboundary and therefore the deformation is infinitesimally trivial. In the last part of this paper we discuss the relation between Mackey's analogue and the quantization problem of the family algebras.
\end{abstract}

\maketitle

\section{Introduction}\label{Introduction}
The classical and quantum \emph{family algebras} are introduced by A. A. Kirillov in the year 2000 in \cite{Ki1} and \cite{Ki2} as a new tool to study representation theory of Lie algebras. They have some interesting properties and various applications in Lie theory, representation theory and mathematical physics.

Let us briefly review the definitions of them. Let  $\mathfrak{g}$ be a finite dimensional complex Lie algebra, $S(\mathfrak{g})$ and $U(\mathfrak{g})$ be the symmetric algebra and the universal enveloping algebra of $\mathfrak{g}$, respectively. Let $G$ be a connected and simply connected Lie group with Lie$(G)=\mathfrak{g}$. $G$ has adjoint actions $\Ad$ on $S(\mathfrak{g})$ and $U(\mathfrak{g})$.

On the other hand, let $V_{\tau}$ be a finite dimensional complex representation of  $\mathfrak{g}$.  Then $\tau$ gives rise to a representation of $G$. Hence $G$ has a natural action on $\End_{\mathbb{C}}V_{\tau}$:
$$
\forall A\in \End_{\mathbb{C}}V_{\tau}, g\in G,~g\cdot A := \tau(g)  A  \tau(g)^{-1}
$$
and therefore $G$ has natural diagonal actions on $\End_{\mathbb{C}}V_{\tau}\otimes_{\mathbb{C}}S(\mathfrak{g})$ and $\End_{\mathbb{C}}V_{\tau}\otimes_{\mathbb{C}}U(\mathfrak{g})$.

Now we come to the definition of the family algebras, see~\cite{Ki1}~and~\cite{Ki2}:
\begin{defi}[The family algebras]\label{defi of family algebras}
The \emph{classical family algebra} is defined to be:
\begin{equation}
\mathcal{C}_{\tau}(\mathfrak{g}):=
(\End_{\mathbb{C}}V_{\tau}\otimes_{\mathbb{C}}S(\mathfrak{g}))^G.
\end{equation}
The \emph{quantum family algebra} is defined to be:
\begin{equation}
\mathcal{Q}_{\tau}(\mathfrak{g}):=
(\End_{\mathbb{C}}V_{\tau}\otimes_{\mathbb{C}}U(\mathfrak{g}))^G.
\end{equation}
\end{defi}

\begin{rmk}
Actually the classical family algebra is a generalization of the invariant polynomial algebra $I(\mathfrak{g})$ and the quantum family algebra is a generalization of the center of the universal enveloping algebra $Z(\mathfrak{g})$, see Example \ref{eg: trivial rep and family algebra} below.
\end{rmk}

Kirillov proves that $\mathcal{C}_{\tau}(\mathfrak{g})$ and $\mathcal{Q}_{\tau}(\mathfrak{g})$ are indeed algebras, that is, they are closed under multiplications. A.A. Kirillov \cite{Ki1}, \cite{Ki2}, A. Joseph \cite{joseph2014modules}, N. Rozhkovskaya \cite{Roz} and M. Tai \cite{tai2013classical} have found various  relations between family algebras and the representations of $\mathfrak{g}$. For other applications, N. Higson relates family algebras with the admissible representations of complex semisimple Lie groups in \cite{Higson2011}

\begin{rmk}
The quantum family algebra is called \emph{the relative Yangian} in \cite{joseph2014modules}.
\end{rmk}

In this paper, we study  the family algebras from another viewpoint. It is well-known that we have a \emph{Poisson bracket} on $S(\mathfrak{g})$ (see \cite{Ki3}): Let $X_i$ be a basis of $\mathfrak{g}$ and $c^k_{ij}$ be the structure constant with respect to the basis $X_i$, then for any $a,b \in S(\mathfrak{g})$, the Poisson bracket is defined to be
$$
\{a,b\}:=c^k_{ij}X_k\cdot\partial^ia\cdot\partial^jb
$$
where $\partial^i$ denotes the partial derivative with respect to $X_i$.

Now we can define the \emph{noncommutative Poisson bracket} on the classical family algebra:
\begin{defi}[The noncommutative Poisson bracket on $\mathcal{C}_{\tau}(\mathfrak{g})$]\label{poisson bracket for cfa}
Let $\mathcal{A}, \mathcal{B}\in \mathcal{C}_{\tau}(\mathfrak{g})$, $\mathcal{A}=A_i\otimes a^i,~\mathcal{B}=B_j\otimes b^j$. We define the noncommutative Poisson bracket $P$ as follows:
\begin{equation}
\{\mathcal{A},\mathcal{B}\}:=A_iB_j\otimes\{a^i, b^j\}.
\end{equation}
\end{defi}

In this paper we will study the properties of the noncommutative Poisson bracket (for short, Poisson bracket) on $\mathcal{C}_{\tau}(\mathfrak{g})$. The following are two important results we get:
\begin{itemize}
\item The Poisson bracket on $\mathcal{C}_{\tau}(\mathfrak{g})$ characterize the first-order $1$-parameter formal deformations from $\mathcal{C}_{\tau}(\mathfrak{g})$ to $\mathcal{Q}_{\tau}(\mathfrak{g})$, just as the Poisson bracket on $S(\mathfrak{g})$ characterize the first-order  $1$-parameter formal deformations from $S(\mathfrak{g})$ to $U(\mathfrak{g})$, see Proposition \ref{first terms of the star product on family algebras}.
\item In the Hochschild cochain complex of $\mathcal{C}_{\tau}(\mathfrak{g})$, the Poisson bracket is a 2-coboundary. In fact we can explicitly find a Hochschild 1-cochain $\nabla$ which maps to the Poisson bracket under the Hochschild differential, see Theorem \ref{poisson bracket is a coboundary}.
\end{itemize}

It is expected that this result has applications in representation theory as well as mathematical physics. For example, it may help us find a quantization map $\mathcal{C}_{\tau}(\mathfrak{g})\rightarrow \mathcal{Q}_{\tau}(\mathfrak{g})$, as proposed by  Higson in \cite{Higson2011}. See Section \ref{section: quantization of family algebras} for details.\\

This paper is organized as follows: In Section \ref{section: fa} we review the family algebras, in Section \ref{section: poisson} we study the first properties of the noncommutative Poisson bracket $P$ on the classical family algebra. In Section \ref{section: deform of fa} we give the relation between $P$ and the  $1$-parameter formal deformation from $\mathcal{C}_{\tau}(\mathfrak{g})$ to $\mathcal{Q}_{\tau}(\mathfrak{g})$. In Section \ref{section: nabla} we prove that the noncommutative Poisson bracket $P$ is a Hochschild $2$-coboundary and therefore the deformation is infinitesimally trivial. In Section \ref{section: quantization of family algebras} we talk about the quantization problem of the family algebras. In the three appendices we summarize the results on Hochschild cohomology, Gerstenhaber bracket and their relation to the deformation theory.

\begin{rmk} Although Kirillov and Higson in \cite{Ki1}, \cite{Ki2} and \cite{Higson2011} require the Lie algebra $\mathfrak{g}$ to be semisimple and the representation $\tau$ to be irreducible, in this paper we do not need this restriction, except for Section \ref{section: quantization of family algebras}.
\end{rmk}

\subsection*{Conventions}
Let $\mathfrak{g}$ be a finite dimensional complex Lie algebra. We use $\{X_i\}$ for a basis for $\mathfrak{g}$ and $c^k_{ij}$ the structure constant of $\mathfrak{g}$ with respect to $\{X_i\}$.

We use the  letters in lower case $a$, $b$ or $a_i$, $b_j$ etc. for the elements in $S(\mathfrak{g})$. The symbol $\partial^k$ denotes the partial derivative of elements in $S(\mathfrak{g})$ with respect to $X_k$.

Let $V_{\tau}$ be a finite representation of $\mathfrak{g}$, we use the capital letter $A$, $B$ or $A_i$, $B_j$ etc. to represent the elements in $\End(V_{\tau})$.

The calligraphic letters $\mathcal{A}$, $\mathcal{B}$ etc. stand for elements in the family algebras $\mathcal{C}_{\tau}(\mathfrak{g})$ and $\mathcal{Q}_{\tau}(\mathfrak{g})$.

In this paper we frequently use the Einstein summation convention: $A_i\otimes a^i$ means taking sum with respect to $i$.

\section{A brief introduction to the Family algebras}\label{section: fa}
We give a quick review of family algebras in this section. Most of the materials in this section can be found in \cite{Ki1} and \cite{Ki2}.\\

First of all, we use the following notation-definition
\begin{defi}\label{tilde family algebras}
\begin{equation}
\begin{split}
\widetilde{\mathcal{C}_{\tau}}(\mathfrak{g})&:=\End_{\mathbb{C}}V_{\tau}\otimes_{\mathbb{C}}S(\mathfrak{g}),\\
\widetilde{\mathcal{Q}_{\tau}}(\mathfrak{g})&:=\End_{\mathbb{C}}V_{\tau}\otimes_{\mathbb{C}}U(\mathfrak{g}).
\end{split}
\end{equation}
\end{defi}

$\widetilde{\mathcal{C}_{\tau}}(\mathfrak{g})$ and $\widetilde{\mathcal{Q}_{\tau}}(\mathfrak{g})$ consist of matrices with entries in $S(\mathfrak{g})$ and $U(\mathfrak{g})$, respectively. Therefore they are algebras in a natural way: for any $A_i\otimes a^i$, $B_j\otimes b^j\in \widetilde{\mathcal{C}_{\tau}}(\mathfrak{g})(\text{ or }\widetilde{\mathcal{Q}_{\tau}}(\mathfrak{g}))$, their product is given by the following formula:
\begin{equation}
(A_i\otimes a^i)\cdot(B_j\otimes b^j):=A_iB_j\otimes a^ib^j.
\end{equation}

\begin{ctn}
$\widetilde{\mathcal{C}_{\tau}}(\mathfrak{g})$ and $\widetilde{\mathcal{Q}_{\tau}}(\mathfrak{g})$  are non-commutative in general.
\end{ctn}

The following simple result will be frequently used:
\begin{lemma}\label{lemma: commute matrices and Lie terms}
In both $\widetilde{\mathcal{C}_{\tau}}(\mathfrak{g})$ and $\widetilde{\mathcal{Q}_{\tau}}(\mathfrak{g})$, the matrix component and the $S(\mathfrak{g})$ component always commute. In more detail, for any $A_i\otimes a^i$, $B_j\otimes b^j\in \widetilde{\mathcal{C}_{\tau}}(\mathfrak{g})(\text{or }\widetilde{\mathcal{Q}_{\tau}}(\mathfrak{g}))$, we have
\begin{equation}
\begin{split}
(A_i\otimes a^i)\cdot (B_j\otimes b^j)=&A_iB_j\otimes a^ib^j\\
=&(\text{Id}\otimes a^i)\cdot (A_iB_j\otimes b^j)\\
=&(A_iB_j\otimes a^i)\cdot(\text{Id}\otimes b^j).
\end{split}
\end{equation}
\end{lemma}
\begin{proof}
It is obvious.
\end{proof}

By Definition \ref{defi of family algebras} we know
$$
\mathcal{C}_{\tau}(\mathfrak{g})=\widetilde{\mathcal{C}_{\tau}}(\mathfrak{g})^G \text{ and } \mathcal{Q}_{\tau}(\mathfrak{g}):=\widetilde{\mathcal{Q}_{\tau}}(\mathfrak{g}))^G.
$$

Now we show that the Lie group action can be reduced to the Lie algebra action.

\begin{prop}[The infinitesimal test for classical family algebra, \cite{Ki1} Section 1]\label{criforcfa}
Let $A_i\otimes a^i\in \widetilde{\mathcal{C}_{\tau}}(\mathfrak{g})$, then $A_i\otimes a^i\in \mathcal{C}_{\tau}(\mathfrak{g})$ if and only if
\begin{equation}\label{criterion for classical family algebra 1}
\forall X\in \mathfrak{g},~[\tau(X), A_i]\otimes a^i+ A_i\otimes \{X,~a^i\}=0.
\end{equation}
In other words,
\begin{equation}\label{criterion for classical family algebra 2}
\forall X\in \mathfrak{g},~[\tau(X), A_i]\otimes a^i=A_i\otimes \{a^i,~X\}.
\end{equation}
\end{prop}
\begin{proof}
By definition \ref{defi of family algebras}, we know that $A_i\otimes a^i\in \mathcal{C}_{\tau}(\mathfrak{g})$ if and only if:
$$
\tau(g)\circ A_i\circ \tau(g)^{-1}\otimes (\ad{g})~ a^i=A_i\otimes a^i.
$$
It is well-known that the adjoint action of $\mathfrak{g}$ on $S(\mathfrak{g})$  is exactly the Poisson bracket. As a result, Equation (\ref{criterion for classical family algebra 1}) and (\ref{criterion for classical family algebra 2}) are infinitesimal versions of the above equation. Since $G$ is connected and simply connected, they are equivalent to the invariance under the action of $G$.
\end{proof}

Similarly we have
\begin{prop}[The infinitesimal test for quantum family algebra, \cite{Ki1} Section 1]\label{criforqfa}
Let $A_i\otimes a^i\in \widetilde{\mathcal{Q}_{\tau}}(\mathfrak{g})$, then $A_i\otimes a^i\in \mathcal{Q}_{\tau}(\mathfrak{g})$  if and only if
$$
\forall X\in \mathfrak{g},~[\tau(X), A_i]\otimes a^i+ A_i\otimes [X,~a^i]=0.
$$
In other words,
$$
\forall X\in \mathfrak{g},~[\tau(X), A_i]\otimes a^i=A_i\otimes [a^i,~X].
$$
\end{prop}
\begin{proof} Similar to the proof of Proposition \ref{criforqfa}. \end{proof}

Then we can prove the following result:

\begin{coro}[see also \cite{Ki1} and \cite{Ki2}]\label{family algebras are subalgebras of tilde}
$\mathcal{C}_{\tau}(\mathfrak{g})$ and $\mathcal{Q}_{\tau}(\mathfrak{g})$ are subalgebras of $\widetilde{\mathcal{C}_{\tau}}(\mathfrak{g})$~and~$\widetilde{\mathcal{Q}_{\tau}}(\mathfrak{g}))$ respectively.
\end{coro}
\begin{proof}
Let $A_i\otimes a^i$ and $B_j\otimes b^j$ be two elements in $\mathcal{C}_{\tau}(\mathfrak{g})$. Their product
$$
(A_i\otimes a^i)\cdot(B_j\otimes b^j)=A_i B_j\otimes a^ib^j.
$$

Now $\forall X\in \mathfrak{g}$,
\begin{align*}
[\tau(X), A_iB_j]\otimes a^ib^j=& [\tau(X), A_i]B_j\otimes a^ib^j+A_i[\tau(X), B_j]\otimes a^ib^j\\
=& ([\tau(X), A_i]\otimes a^i)\cdot (B_j\otimes b^j)+(A_i\otimes a^i)\cdot([\tau(X), B_j]\otimes b^j).
\end{align*}
The second equality is because of Lemma \ref{lemma: commute matrices and Lie terms}: the matrix component always \textbf{commutes} with the $S(\mathfrak{g})$ component.

Now by Proposition $\ref{criforcfa}$, we know
\begin{align*}
\text{the above formula }=&(A_i\otimes \{a^i, X\})\cdot (B_j\otimes b^j)+(A_i\otimes a^i)\cdot B_j\otimes \{b^j, X\}\\
=&A_iB_j\otimes \{a^i, X\}b^j+A_iB_j\otimes a^i\{b^j, X\}~(\text{Proposition}~\ref{criforcfa})\\
=&A_iB_j\otimes \{a^ib^j, X\}.
\end{align*}

Hence we get
$$A_i B_j\otimes a^ib^j\in \mathcal{C}_{\tau}(\mathfrak{g}).$$

In the same way we can show that if $A_i\otimes a^i$ and $B_j\otimes b^j$ are in $\mathcal{Q}_{\tau}(\mathfrak{g})$,
then
$$
A_i B_j\otimes a^ib^j\in \mathcal{Q}_{\tau}(\mathfrak{g}).
$$
\end{proof}

It is not difficult to see that the family algebras contains non-zero elements. In fact, let $I(\mathfrak{g})=S(\mathfrak{g})^{\mathfrak{g}}$ be the invariant subalgebra of $S(\mathfrak{g})$ and $Z(\mathfrak{g})$ be the center of $U(\mathfrak{g})$. We have

\begin{prop}[\cite{Ki2}]
$I(\mathfrak{g})$ embeds into $\mathcal{C}_{\tau}(\mathfrak{g})$ as scalar matrices
\begin{equation}
\begin{split}
I(\mathfrak{g})&\hookrightarrow\mathcal{C}_{\tau}(\mathfrak{g})\\
a &\mapsto \text{Id}\otimes a.
\end{split}
\end{equation}
Similarly $Z(\mathfrak{g})$ embeds into $\mathcal{Q}_{\tau}(\mathfrak{g})$ as scalar matrices too.
\end{prop}
\begin{proof} It is obvious that $I(\mathfrak{g})$ embeds into $\widetilde{\mathcal{C}_{\tau}}(\mathfrak{g})$ as scalar matrices. Now by Proposition \ref{criforcfa}, it is easy to see that the image is contained in $\mathcal{C}_{\tau}(\mathfrak{g})$.

The proof for $Z(\mathfrak{g})$ and $\mathcal{Q}_{\tau}(\mathfrak{g})$ is the same. \end{proof}

\begin{eg}\label{eg: trivial rep and family algebra}
For any $\mathfrak{g}$, when the representation $\tau$ is the trivial representation, we see that $I(\mathfrak{g})=\mathcal{C}_{\tau}(\mathfrak{g})$ and $Z(\mathfrak{g})=\mathcal{Q}_{\tau}(\mathfrak{g})$.
\end{eg}

\begin{eg}\label{sl2 standard rep}
For $\mathfrak{g}=\sll(2,\mathbb{C})$  and $\{e,f,h\}$ be the standard basis of $\sll(2,\mathbb{C})$ which satisfies the commutation relation
\begin{equation}
[e,f]=h,~~[h,e]=2e,~~[h,f]=-2f.
\end{equation}
Let $\tau$ be the $2$-dimensional standard representation, we can find an element $M\in \mathcal{C}_{\tau}(\mathfrak{g})$ which is not in $I(\mathfrak{g})$. In fact
\begin{equation}
M=\left.\begin{pmatrix}\frac{h}{2} & f \\ e & -\frac{h}{2}\end{pmatrix}\right.
\end{equation}
We can also find an element in $\mathcal{Q}_{\tau}(\mathfrak{g})$ with the same expression of $M$, see \cite{Ki1} and \cite{Ki2}.
\end{eg}

\begin{rmk} When $\tau$ is nontrivial irreducible and $\mathfrak{g}$ is semisimple, $I(\mathfrak{g})$ is not equal to  $\mathcal{C}_{\tau}(\mathfrak{g})$ and $Z(\mathfrak{g})$ is not equal to $\mathcal{Q}_{\tau}(\mathfrak{g})$ either, see Corollary \ref{famili algebra not only scalar} below or \cite{Ki1}.
\end{rmk}

\section{The noncommutative Poisson bracket on  $\widetilde{\mathcal{C}_{\tau}}(\mathfrak{g})$}\label{section: poisson}

The noncommutative Poisson bracket on $\mathcal{C}_{\tau}(\mathfrak{g})$ in Definition \ref{poisson bracket for cfa} can be automatically extended to $\widetilde{\mathcal{C}_{\tau}}(\mathfrak{g})$:
\begin{defi}\label{poisson bracket for tilde cfa}
Let $\mathcal{A}, \mathcal{B}\in \widetilde{\mathcal{C}_{\tau}}(\mathfrak{g})$, $\mathcal{A}=A_i\otimes a^i,~\mathcal{B}=B_j\otimes b^j$. We define the noncommutative Poisson bracket as follows:
\begin{equation}
\{\mathcal{A},\mathcal{B}\}:=A_iB_j\otimes\{a^i, b^j\}.
\end{equation}

We will also denote the noncommutative Poisson bracket by $P$.
\end{defi}

\begin{rmk}
P. Xu gives a similar construction in \cite{PingXu1994} Example 2.2.
\end{rmk}

\begin{ctn}
The Poisson bracket on $\widetilde{\mathcal{C}_{\tau}}(\mathfrak{g})$ is not anti-symmetric and does not satisfies the Leibniz rule and the Jacobi identity.
\end{ctn}

Nevertheless, J. Block and E. Getzler in 1992 give a definition of Poisson brackets on noncommutative algebras in \cite{BlockGetzler1992} and we can prove that our noncommutative Poisson bracket $P$ satisfies the requirement of Poisson bracket in that sense:

\begin{defi}[\cite{BlockGetzler1992} Definition 1.1]\label{defi of noncomm Poisson bracket in general}
A \emph{Poisson bracket} on a (possibly noncommutative) algebra $A$ is a Hochschild $2$-cocycle $P\in Z^2(A,A)$ such that $P\circ P \in C^3(A,A)$ is a $3$-coboundary. In other words
\begin{equation}
P\circ P \in B^3(A,A)\subset Z^3(A,A)\subset C^3(A,A).
\end{equation}
\end{defi}

For Hochschild cohomology see Appendix \ref{app:Hochschild cohomology} and for the definition of $P\circ P$ see Proposition \ref{circ on 2-cochain}.

\begin{rmk}
In Definition \ref{defi of noncomm Poisson bracket in general}, we may take the condition $P\in Z^2(A,A)$  as a noncommutative Leibniz rule and $P\circ P \in B^3(A,A)$ as a noncommutative Jacobi identity. They together implies that $P$ can be lift to an associative product on $A$ \emph{up to order $3$}, see Corollary \ref{deformation up to 2 and 3}.
\end{rmk}

For our algebra $\widetilde{\mathcal{C}_{\tau}}(\mathfrak{g})$ and the Poisson bracket $P$ in Definition \ref{poisson bracket for tilde cfa}, first we can prove that $P$ is a $2$-cocycle. We have the following proposition:

\begin{prop}\label{generalized Leibniz rule}
For any $\mathcal{A}, \mathcal{B}, \mathcal{C} \in\widetilde{\mathcal{C}_{\tau}}(\mathfrak{g})$, we have
\begin{equation}
\mathcal{A\{B,C\}-\{AB,C\}+\{A,BC\}-\{A,B\}C}=0.
\end{equation}

In other words, we have $\dH P=0$ where $\dH$ is the Hochschild differential operator. Therefore
\begin{equation}
P\in Z^2(\widetilde{\mathcal{C}_{\tau}}(\mathfrak{g}),\widetilde{\mathcal{C}_{\tau}}(\mathfrak{g})).
\end{equation}
\end{prop}
\begin{proof}  Let $\mathcal{A}=A_i\otimes a^i$, $\mathcal{B}=B_j\otimes b^j$ and $\mathcal{C}=C_k\otimes c^k$. Recall that we have Lemma \ref{lemma: commute matrices and Lie terms}: the matrix component and the $S(\mathfrak{g})$ component always commute. Then by the definition of the Poisson bracket we get
\begin{align*}
&\mathcal{A\{B,C\}-\{AB,C\}+\{A,BC\}-\{A,B\}C}\\
=& A_iB_jC_k\otimes (a_i\{b_j,c_k\}-\{a_ib_j,c_k\}+\{a_i,b_jc_k\}-\{a_i,b_j\}c_k).
\end{align*}
By the Leibniz rule of the (ordinary) Poisson bracket on $S(\mathfrak{g})$ we know that
$$
a_i\{b_j,c_k\}-\{a_ib_j,c_k\}+\{a_i,b_jc_k\}-\{a_i,b_j\}c_k=0.
$$
\end{proof}

\begin{prop}\label{generalized Jacobi id}
$P\circ P$ is a $3$-coboundary. In other words, $P\circ P\in B^3(\widetilde{\mathcal{C}_{\tau}}(\mathfrak{g}),\widetilde{\mathcal{C}_{\tau}}(\mathfrak{g}))$.
\end{prop}
\begin{proof}
In fact, we can define a $2$-cochain $\Phi\in C^2(\widetilde{\mathcal{C}_{\tau}}(\mathfrak{g}),\widetilde{\mathcal{C}_{\tau}}(\mathfrak{g}))$ as follows: Let $\mathcal{A}=A\otimes a$ and $\mathcal{B}=B\otimes b$ (to simplify the notation we omit the super and sub-indices)
\begin{equation}
\begin{split}
\Phi(\mathcal{A},\mathcal{B}):=&AB\otimes\frac{1}{2}c^s_{ij}c^t_{kl}X_s\cdot X_t\cdot \partial^i\partial^ka\cdot \partial^j\partial^lb\\
+&AB\otimes\frac{1}{3}c^t_{ks}c^s_{ji}X_t(\partial^k\partial^ja\cdot \partial^ib+\partial^ia\cdot \partial^k\partial^jb).
\end{split}
\end{equation}
Then we have
\begin{equation}
P\circ P+\dH \Phi=0.
\end{equation}
 For any $\mathcal{A}, \mathcal{B}, \mathcal{C} \in\widetilde{\mathcal{C}_{\tau}}(\mathfrak{g})$, by Definition
\begin{equation}
P\circ P(\mathcal{A}, \mathcal{B}, \mathcal{C})=\{\mathcal{A},\{\mathcal{B},\mathcal{C}\}\}-\{\{\mathcal{A},\mathcal{B}\},\mathcal{C}\}.
\end{equation}

Let $\mathcal{A}=A\otimes a$, $\mathcal{B}=B\otimes b$ and $\mathcal{C}=C\otimes c$ , then
$$
P\circ P(\mathcal{A}, \mathcal{B}, \mathcal{C})= ABC\otimes (\{a,\{b,c\}\}-\{\{a,b\},c\})\\.
$$

Now the problem reduces to $S(\mathfrak{g})$. We have the following lemma:
\begin{lemma}\label{lemma for HH3 cobdy in S(g)}
We can define a $2$-cochain $\phi\in C^2(S(\mathfrak{g}),S(\mathfrak{g}))$ as follows: for any $a, b \in S(\mathfrak{g})$
\begin{equation}
\phi(a,b):=\frac{1}{2}c^s_{ij}c^t_{kl}X_s\cdot X_t\cdot \partial^i\partial^ka \cdot \partial^j\partial^lb+\frac{1}{3}c^t_{ks}c^s_{ji}X_t\cdot (\partial^k\partial^ja\cdot \partial^ib+\partial^ia\cdot \partial^k\partial^jb).
\end{equation}

Then for any $a, b, c \in S(\mathfrak{g})$ we have
\begin{equation}
\{a,\{b,c\}\}-\{\{a,b\},c\}+(\dH\, \phi)(a,b,c)=0.
\end{equation}
\end{lemma}
\begin{proof}[Proof of Lemma \ref{lemma for HH3 cobdy in S(g)}] We can check it by hand using Jacobi identity. Another approach involves the star-porduct on $S(\mathfrak{g})$ and the general result of deformation theory and we defer it to Proposition \ref{prop for HH3 cobdy in S(g)}. \end{proof}

Now we have Lemma \ref{lemma for HH3 cobdy in S(g)}. By abusing the notations we have $\Phi=\text{Id}\otimes \phi$, then we immediately get $P\circ P+\dH \Phi=0$. This finishes the proof of Proposition \ref{generalized Jacobi id}. \end{proof}

By Proposition \ref{generalized Leibniz rule} and \ref{generalized Jacobi id} we know that the noncommutative Poisson bracket in Definition \ref{poisson bracket for tilde cfa} is indeed a Poisson bracket in the sense of Definition \ref{defi of noncomm Poisson bracket in general}.

Before we move on, we need to prove that the Poisson bracket indeed maps $\mathcal{C}_{\tau}(\mathfrak{g})\otimes \mathcal{C}_{\tau}(\mathfrak{g})$ to $\mathcal{C}_{\tau}(\mathfrak{g})$. That is the following proposition:

\begin{prop}\label{poisson bracket mas cfa to cfa}
For any $\mathcal{A}, \mathcal{B}\in \mathcal{C}_{\tau}(\mathfrak{g})$, we have that $\{\mathcal{A}, \mathcal{B}\}$ is still in $\mathcal{C}_{\tau}(\mathfrak{g})$. In other words, the noncommutative Poisson bracket in Definition \ref{poisson bracket for cfa} is well-defined.
\end{prop}
\begin{proof} We can proof this proposition by computation using Proposition \ref{criforcfa} and the definition of the noncommutative Poisson bracket $P$. In Section \ref{section: nabla} we will give another proof using a different construction of $P$. See Corollary \ref{coro: poisson bracket mas cfa to cfa}. \end{proof}

\section{The noncommutative Poisson bracket and the  $1$-parameter formal  deformation of $\widetilde{\mathcal{C}_{\tau}}(\mathfrak{g})$}\label{section: deform of fa}
In this section we will show that the Poisson bracket plays an essential role in the  $1$-parameter formal  deformation of $\widetilde{\mathcal{C}_{\tau}}(\mathfrak{g})$.

\subsection{A quick review of the  $1$-parameter formal  deformation from $S(\mathfrak{g})$ to $U(\mathfrak{g})$ and the Poisson bracket}

Before studying the deformation of $\widetilde{\mathcal{C}_{\tau}}(\mathfrak{g})$, let us first review the corresponding theory of $S(\mathfrak{g})$ and  $U(\mathfrak{g})$.

It is well-known that
\begin{equation}
S(\mathfrak{g})=\text{T}(\mathfrak{g})/(X\otimes Y-Y\otimes X)
\end{equation}
and
\begin{equation}
U(\mathfrak{g})=\text{T}(\mathfrak{g})/(X\otimes Y-Y\otimes X-[X,Y])
\end{equation}
where $\text{T}(\mathfrak{g})$ is the tensor algebra of $\mathfrak{g}$.

Moreover, we consider the algebra
\begin{equation}
U_t(\mathfrak{g})=\text{T}(\mathfrak{g})/(X\otimes Y-Y\otimes X-t[X,Y]).
\end{equation}
For $t\neq 0$ all the algebras $U_t(\mathfrak{g})$ are isomorphic to $U(\mathfrak{g})$, and when $t=0$, $U_0(\mathfrak{g})$ is isomorphic to $S(\mathfrak{g})$. $t$ is called the \emph{deformation parameter}.

\begin{rmk}
For more details on the $1$-parameter formal deformation of an associative algebra see Appendix \ref{app:HH and deformation}.
\end{rmk}

We have the \emph{Poincar\'{e}-Birkhoff-Witt} map ($\IPBW$) from $S(\mathfrak{g})$ to $U_t(\mathfrak{g})$ given by:
\begin{equation}
\begin{split}
\IPBW:~~~~~~S(\mathfrak{g}) ~~~~~~&\longrightarrow ~~~~~~U_t(\mathfrak{g})\\
X_1X_2\ldots X_k &\longmapsto \sum_{\sigma\in S_k}\frac{1}{k!}X_{\sigma(1)}X_{\sigma(2)}\ldots X_{\sigma(k)}.
\end{split}
\end{equation}
The Poincar\'{e}-Birkhoff-Witt theorem (see \cite{Knapp2002}) tells us that the above map $\IPBW$ is an isomorphism between \textbf{$\mathfrak{g}$-vector spaces}.

\begin{rmk}
The map $\IPBW$ is not an algebraic isomorphism unless $\mathfrak{g}$ is an abelian Lie algebra or $t=0$.
\end{rmk}

Therefore we have the following definition

\begin{defi}\label{star product on S(g)}
The map $\IPBW$ pulls back the multiplication of $U_t(\mathfrak{g})$ to $S(\mathfrak{g})$ and we call it the \textbf{star-product} on $S(\mathfrak{g})$, denoted by $*_t$. For any $a,b\in S(\mathfrak{g})$
\begin{equation}
a *_t b:= \IPBW^{-1}\,(\IPBW(a)\cdot \IPBW(b)).
\end{equation}
In particular we denote $*_1$ simply as $*$. When $t=0$ the star-product reduces to the original production on $S(\mathfrak{g})$.
Obviously $*_t$ satisfies the associativity law because the multiplication on $U_t(\mathfrak{g})$ is associative.
\end{defi}

Now by definition, the map $\IPBW$ gives an \emph{algebraic isomorphism}
\begin{equation}
\IPBW: (S(\mathfrak{g}), *_t)\overset{\sim}\longrightarrow (U_t(\mathfrak{g}), \cdot).
\end{equation}
Therefore we can identify $U_t(\mathfrak{g})$ with $(S(\mathfrak{g}), *_t)$, especially we can identify $U(\mathfrak{g})$ with $(S(\mathfrak{g}), *)$.

\begin{rmk}
Our star-product $*_t$ is not exactly the same as the sart-product constructed by Kontsevich in \cite{Kontsevich2003} Section 8. Nevertheless, they give isomorphic algebra structures on $S(\mathfrak{g})$.
\end{rmk}

The star-product $*_t$ depends on the deformation parameter $t$. In fact we can write the first few terms of $*_t$ .
\begin{prop}[\cite{Gutt1983} Section 3]\label{first terms of the star product on S(g)}
We can write $*_t$ as
\begin{equation}
a *_t b= ab+\frac{t}{2}\{a,b\}+O(t^2).
\end{equation}
\end{prop}
\begin{proof}
This result is well-known and can be found in, say, \cite{Gutt1983} Section 3.
\end{proof}

 In other words, the Poisson bracket on $S(\mathfrak{g})$ is exactly the \emph{first-order}  $1$-parameter formal deformation from $S(\mathfrak{g})$ to $U(\mathfrak{g})$.

 \begin{rmk}
In fact we can find the expressing of the $t^2$ term in the star-product. According to \cite{Grabowski1992} Remark 4.7, for any $a, b\in S(\mathfrak{g})$, the $t^2$ term is
\begin{equation}\label{m2 in star}
m_2(a,b):=\frac{1}{8}c^s_{ij}c^t_{kl}X_s\cdot X_t\cdot \partial^i\partial^ka \cdot \partial^j\partial^lb+\frac{1}{12}c^t_{ks}c^s_{ji}X_t\cdot (\partial^k\partial^ja\cdot \partial^ib+\partial^ia\cdot \partial^k\partial^jb).
\end{equation}
\end{rmk}

Now we can give another proof of Lemma \ref{lemma for HH3 cobdy in S(g)}
\begin{prop}[Lemma \ref{lemma for HH3 cobdy in S(g)}]\label{prop for HH3 cobdy in S(g)}
We can define a $2$-cochain $\phi\in C^2(S(\mathfrak{g}),S(\mathfrak{g}))$ as follows: for any $a, b \in S(\mathfrak{g})$
\begin{equation}
\phi(a,b):=\frac{1}{2}c^s_{ij}c^t_{kl}X_s\cdot X_t\cdot \partial^i\partial^ka \cdot \partial^j\partial^lb+\frac{1}{3}c^t_{ks}c^s_{ji}X_t\cdot (\partial^k\partial^ja\cdot \partial^ib+\partial^ia\cdot \partial^k\partial^jb).
\end{equation}

Then for any $a, b, c \in S(\mathfrak{g})$ we have
\begin{equation}\label{P phi compatible}
\{a,\{b,c\}\}-\{\{a,b\},c\}+(\dH\, \phi)(a,b,c)=0.
\end{equation}
\end{prop}
\begin{proof} The result is clear in the framework of deformation theory (see Appendix \ref{app:HH and deformation}). Let $m=*_t$ be the star-product. Compare Proposition (\ref{first terms of the star product on S(g)}), Equation (\ref{m2 in star}) and Equation (\ref{deformation of algebra2}) we get
$$
P=2m_1 \text{~~and~~} \phi=4m_2.
$$
where we denote the Poisson bracket on $S(\mathfrak{g})$ by $P$ too.

Since we know from the definition that the star-product is associative, by Proposition \ref{deformation up to 2 and 3} and in particular Equation (\ref{dH m 2 and circ}) we get
$$
m_1\circ m_1+\dH m_2=0
$$
hence
$$
P\circ P+ \dH \phi=0
$$
and this is exactly Equation (\ref{P phi compatible}).
\end{proof}

If we restrict ourselves to the invariant subalgebra $I(\mathfrak{g})=S(\mathfrak{g})^{\mathfrak{g}}$, then we have the following well-known result about the Poisson bracket:

\begin{prop}[\cite{Ki3}]\label{vanish of Poisson bracket on I(g)}
The Poisson bracket vanishes on $I(\mathfrak{g})$. In other words, for any $a,b \in I(\mathfrak{g})$, we have
\begin{equation}
\{a,b\}=0.
\end{equation}
\end{prop}
\begin{proof} This result comes from the definition of $I(\mathfrak{g})$. In fact
\begin{align*}
\{a,b\} &= c^k_{ij}X_k (\partial^i a)(\partial^j b)\\
&=(c^k_{ij}X_k \partial^i a)(\partial^j b)\\
&=(\ad X_j\,(a)) (\partial^j b).
\end{align*}

Since $a \in I(\mathfrak{g})$, we know that $\ad X_j\,(a)=0$ for any $X_j$, as a result, $\{a,b\}=0$.
\end{proof}

On the other hand, we can easily find the image of $\IPBW$ restricted on $I(\mathfrak{g})$.

\begin{prop}[\cite{Knapp2002}]\label{Ipbw I(g) to Z(g)}
The image of $I(\mathfrak{g})$ under the Poincar\'{e}-Birkhoff-Witt map  $\IPBW$ is exactly $Z_t(\mathfrak{g})$, the center of $U_t(\mathfrak{g})$. In other words,
\begin{equation}
\IPBW: I(\mathfrak{g})\rightarrow Z_t(\mathfrak{g})
\end{equation}
is an isomorphism between vector spaces.
\end{prop}
\begin{proof}
Keep in mind that $\IPBW:~ S(\mathfrak{g})\rightarrow U_t(\mathfrak{g})$ is an isomorphism between $\mathfrak{g}$-vector spaces,i.e. it is compatible with the $\mathfrak{g}$-actions. \end{proof}

\begin{rmk}\label{vanish of 1st order deformation on I(g)}
Proposition \ref{vanish of Poisson bracket on I(g)} and Proposition \ref{Ipbw I(g) to Z(g)} tell us that the first-order  $1$-parameter formal deformation from $I(\mathfrak{g})$ to $Z(\mathfrak{g})$ is zero.
\end{rmk}

In fact we have the much deeper \emph{Duflo's isomorphism theorem}:

\begin{thm}[\cite{Duflo1977}, \cite{Kontsevich2003} Section 8, \cite{AleMein2000} and \cite{CalRoss2011}]\label{Duflo's isomorphism}
There exists an \textbf{algebraic isomorphism}:
\begin{equation}
\text{Duf}: I(\mathfrak{g})\rightarrow Z(\mathfrak{g})
\end{equation}
\end{thm}
We do not give the proof here. Interested readers can check the references above.

\begin{rmk}
In general, the map Duf will be different from the Poincar\'{e}-Birkhoff-Witt map  $\IPBW$ in Proposition \ref{Ipbw I(g) to Z(g)}, although they have the same domain and image.
\end{rmk}

\begin{rmk}
Since $Z(\mathfrak{g})$ is isomorphic to $Z_t(\mathfrak{g})$ as algebras, the map Duf can be easily generalized to the map $\text{Duf}_t: I(\mathfrak{g})\rightarrow Z_t(\mathfrak{g})$ for any $t$.
\end{rmk}

\subsection{The  $1$-parameter formal deformation from $\widetilde{\mathcal{C}_{\tau}}(\mathfrak{g})$ to $\widetilde{\mathcal{Q}_{\tau}}(\mathfrak{g})$  and the noncommutative Poisson bracket}

We can generalize the constructions of $S(\mathfrak{g})$ to $\widetilde{\mathcal{C}_{\tau}}(\mathfrak{g})$ in this subsection.

\begin{defi}\label{defi of deformed family algebras}
We define the algebra $\widetilde{\mathcal{Q}^t_{\tau}}(\mathfrak{g})$ as
\begin{equation}
\widetilde{\mathcal{Q}^t_{\tau}}(\mathfrak{g}):=\End V_{\tau}\otimes U_t(\mathfrak{g}).
\end{equation}

Moreover, we define
\begin{equation}
\mathcal{Q}^t_{\tau}(\mathfrak{g}):=(\End V_{\tau}\otimes U_t(\mathfrak{g}))^G.
\end{equation}
\end{defi}

By definition, we have $\widetilde{\mathcal{Q}^0_{\tau}}(\mathfrak{g})=\widetilde{\mathcal{C}_{\tau}}(\mathfrak{g})$, $\mathcal{Q}^0_{\tau}(\mathfrak{g})=\mathcal{C}_{\tau}(\mathfrak{g})$ and for any $t\neq 0$ we have $\widetilde{\mathcal{Q}^t_{\tau}}(\mathfrak{g})\cong\widetilde{\mathcal{Q}_{\tau}}(\mathfrak{g})$, $\mathcal{Q}^t_{\tau}(\mathfrak{g})\cong \mathcal{Q}_{\tau}(\mathfrak{g})$.

We also have the Poincar\'{e}-Birkhoff-Witt map on the family algebras:

\begin{defi}\label{defi of PBW on family algebras}
The Poincar\'{e}-Birkhoff-Witt map $\FPBW$ on family algebras is defined to be $\text{Id} \otimes\IPBW$. In other words:
\begin{equation}
\begin{split}
\FPBW: ~~~\widetilde{\mathcal{C}_{\tau}}(\mathfrak{g})&\longrightarrow \widetilde{\mathcal{Q}^t_{\tau}}(\mathfrak{g})\\
A_i\otimes a^i &\longmapsto A_i\otimes \IPBW(a^i).
\end{split}
\end{equation}

$\FPBW$ is an isomorphism between $\mathfrak{g}$-vector spaces.
\end{defi}

As $\IPBW$, $\FPBW$ is not an algebraic isomorphism either. Nevertheless it can also pull back the product on $\widetilde{\mathcal{Q}^t_{\tau}}(\mathfrak{g})$ to $\widetilde{\mathcal{C}_{\tau}}(\mathfrak{g})$:

\begin{defi}\label{defi of star product on family algebras}
The star-product $*_t$ on $\widetilde{\mathcal{C}_{\tau}}(\mathfrak{g})$ is defined to the pull-back of the product on $\widetilde{\mathcal{Q}^t_{\tau}}(\mathfrak{g})$ via the map $\FPBW$. In other words, for any $\mathcal{A}, \mathcal{B}\in \widetilde{\mathcal{C}_{\tau}}(\mathfrak{g})$
\begin{equation}
\mathcal{A} *_t \mathcal{B}:=\FPBW^{-1}\,(\FPBW(\mathcal{A})\cdot \FPBW(\mathcal{B})).
\end{equation}
Moreover, if we write $\mathcal{A}=A_i\otimes a^i$ and $\mathcal{B}=B_j\otimes b^j$, then
\begin{equation}\label{star product components}
(A_i\otimes a^i )*_t (B_j\otimes b^j)=A_iB_j\otimes (a^i*_t b^j).
\end{equation}
\end{defi}

Therefore the map $\FPBW$ gives an \emph{algebraic isomorphism}
\begin{equation}
\FPBW: (\widetilde{\mathcal{C}_{\tau}}(\mathfrak{g}), *_t)\overset{\sim}\longrightarrow (\widetilde{\mathcal{Q}^t_{\tau}}(\mathfrak{g}), \cdot).
\end{equation}
Hence we can identify $\widetilde{\mathcal{Q}^t_{\tau}}(\mathfrak{g})$ with $(\widetilde{\mathcal{C}_{\tau}}(\mathfrak{g}), *_t)$, especially we can identify $\widetilde{\mathcal{Q}_{\tau}}(\mathfrak{g})$ with $(\widetilde{\mathcal{C}_{\tau}}(\mathfrak{g}), *)$.

For the star-product on $\widetilde{\mathcal{C}_{\tau}}(\mathfrak{g})$, we also have

\begin{prop}\label{first terms of the star product on family algebras}
We can write the star-product $*_t$ on $\widetilde{\mathcal{C}_{\tau}}(\mathfrak{g})$ as
\begin{equation}
\mathcal{A} *_t \mathcal{B}= \mathcal{A} \mathcal{B}+\frac{t}{2}\{\mathcal{A}, \mathcal{B}\}+O(t^2).
\end{equation}

In other words, the Poisson bracket on $\widetilde{\mathcal{C}_{\tau}}(\mathfrak{g})$ is exactly the first-order $1$-parameter formal  deformation from $\widetilde{\mathcal{C}_{\tau}}(\mathfrak{g})$ to$\widetilde{\mathcal{Q}_{\tau}}(\mathfrak{g})$.
\end{prop}
\begin{proof} This is just a combination of the definition of star-product (Definition \ref{defi of star product on family algebras}), the definition of noncommutative Poisson bracket (Definition \ref{poisson bracket for tilde cfa}) and Proposition \ref{first terms of the star product on S(g)}. \end{proof}
\begin{rmk}
By now, the results in this subsection exemplified the slogan "the deformation theory of an algebra $A$ is the same as that of the matrix algebra $\text{Mat}_{n\times n}(A)$." However, when restrict to the invariant subalgebras, these two become different.
\end{rmk}

If we restrict ourselves to the family algebra $\mathcal{C}_{\tau}(\mathfrak{g})$, i.e. the invariant subalgebra of $\widetilde{\mathcal{C}_{\tau}}(\mathfrak{g})$, we get the follow proposition which is similar to Proposition \ref{Ipbw I(g) to Z(g)}

\begin{prop}\label{Fpbw C(g) to Q(g)}
The image of $\mathcal{C}_{\tau}(\mathfrak{g})$ under the Poincar\'{e}-Birkhoff-Witt map  $\FPBW$ is exactly $\mathcal{Q}^t_{\tau}(\mathfrak{g})$, the invariant subalgebra of $\widetilde{\mathcal{Q}^t_{\tau}}(\mathfrak{g})$. In other words,
\begin{equation}
\FPBW: \mathcal{C}_{\tau}(\mathfrak{g})\rightarrow \mathcal{Q}^t_{\tau}(\mathfrak{g})
\end{equation}
is an isomorphism between vector spaces.
\end{prop}
\begin{proof} Just remember that $\FPBW:~ \widetilde{\mathcal{C}_{\tau}}(\mathfrak{g})\rightarrow \widetilde{\mathcal{Q}^t_{\tau}}(\mathfrak{g})$ is an isomorphism between $\mathfrak{g}$-vector spaces,i.e. it is compatible with the $\mathfrak{g}$-actions.\end{proof}

Now it is natural to ask for the corresponding result of Proposition \ref{vanish of Poisson bracket on I(g)} and the Duflo's isomorphism theorem \ref{Duflo's isomorphism} on family algebras.

In fact, in Theorem \ref{poisson bracket is a coboundary} of this paper we will prove that the noncommutative Poisson bracket vanishes in the \emph{Hochschild cohomology}. The generalization of Duflo's isomorphism theorem to family algebras is still an open problem, see Section \ref{section: quantization of family algebras}.

\section{The vanishing of the noncommutative Poisson bracket in $\HH^2(\mathcal{C}_{\tau}(\mathfrak{g}))$}\label{section: nabla}

\subsection{The twisted gradient map}
In this section we focus on the classical family algebra $\mathcal{C}_{\tau}(\mathfrak{g})$ and the matrix algebra $\widetilde{\mathcal{C}_{\tau}}(\mathfrak{g})$.

\begin{defi} [The twisted gradient map]\label{defofnabla}
We define a map $\nabla:~\widetilde{\mathcal{C}_{\tau}}(\mathfrak{g})\rightarrow \widetilde{\mathcal{C}_{\tau}}(\mathfrak{g})$~as follows: Fix a basis  $X_k$ of $\mathfrak{g}$. Let $\mathcal{A}=A_i\otimes a^i\in \widetilde{\mathcal{C}_{\tau}}(\mathfrak{g})$, then
\begin{equation}
\nabla(A_i\otimes a^i):=A_i\tau(X_k)\otimes \partial^k(a^i).
\end{equation}

Notice that $\widetilde{\mathcal{C}_{\tau}}(\mathfrak{g})$ is nothing but a matrix algebra with entries in $S(\mathfrak{g})$. In the form of matrices,
\begin{equation}\label{nabla as matrices}
\nabla(\mathcal{A})=\partial^k(\mathcal{A})\tau(X_k).
\end{equation}
Hence $\nabla$ is a first-order differential operator on $\widetilde{\mathcal{C}_{\tau}}(\mathfrak{g})$ and we call it the \emph{twisted gradient map}.
\end{defi}

From Equation (\ref{nabla as matrices}) it is not difficult to see that the map $\nabla$ dose not depend on the concrete expression of $\mathcal{A}\in \widetilde{\mathcal{C}_{\tau}}(\mathfrak{g})$ as $A_i\otimes a^i$.

To show $\nabla$ is a well-defined map, it is now sufficient to prove the following proposition:
\begin{prop}\label{nabla indep of basis}
The map $\nabla: \widetilde{\mathcal{C}_{\tau}}(\mathfrak{g})\rightarrow \widetilde{\mathcal{C}_{\tau}}(\mathfrak{g})$ is independent of the choice of the basis of $\mathfrak{g}$.
\end{prop}
\begin{proof} We need to do some computations. Let  $\tilde{X}_j$  be another basis of  $\mathfrak{g}$. Then
$$
\tilde{X}_j=T^k_j X_k
$$
where $T^k_j$ is the transition matrix. Then, let $\tilde{\partial}^j$ be the partial derivation with respect to $\tilde{X}_j$, we have
$$
\tilde{\partial}^j=(T^{-1})^j_k \partial^k.
$$
Let $\tilde{\nabla}$ be the $\nabla$ map under the basis $\tilde{X}_j$, for $A_i\otimes a^i\in\widetilde{\mathcal{C}_{\tau}}(\mathfrak{g})$, we have
\begin{align*}
\tilde{\nabla}(A_i\otimes a^i)=&A_i\tau(\tilde{X}_j)\otimes \tilde{\partial}^j(a^i)\\
=& A_i\tau(T^k_j X_k)\otimes (T^{-1})^j_l \partial^l(a^i).
\end{align*}
The constant $(T^{-1})^j_l$ can be moved to the first component, hence
\begin{align*}
\text{the above}=& T^k_j (T^{-1})^j_l A_i\tau(X_k)\otimes \partial^l(a^i)\\
=& \delta^k_lA_i\tau(X_k)\otimes \partial^l(a^i)\\
=& A_i\tau(X_k)\otimes \partial^k(a^i)\\
=& \nabla(A_i\otimes a^i)
\end{align*}
So $\nabla$ is invariant under the change of basis of $\mathfrak{g}$.
\end{proof}

The map $\nabla$ is obviously $\mathbb{C}$-linear, moreover it has the following important property:

\begin{prop}\label{nabla maps cfa to cfa}
The image  under $\nabla$ of the subalgebra $\mathcal{C}_{\tau}(\mathfrak{g})$ is contained in $\mathcal{C}_{\tau}(\mathfrak{g})$ itself.
\end{prop}
\begin{proof}
The proof requires some careful computations.

Let $A_i\otimes a^i\in \mathcal{C}_{\tau}(\mathfrak{g})$, then
$$
\nabla(A_i\otimes a^i)=A_i\tau(X_k)\otimes \partial^k(a^i).
$$

By the infinitesimal test of the classical family algebra as in Proposition \ref{criforcfa}, it is sufficient to show that  for $X_j$ which is one of the basis of $\mathfrak{g}$, we have
\begin{equation}\label{eq of image}
[\tau(X_j), A_i\tau(X_k)]\otimes \partial^ka^i=A_i\tau(X_k)\otimes \{\partial^ka^i, X_j\}.
\end{equation}

In fact
\begin{align*}
\text{the left hand side of Equation (\ref{eq of image})}\\=&[\tau(X_j), A_i]\tau(X_k)\otimes \partial^ka^i+A_i[\tau(X_j),\tau(X_k)]\otimes \partial^ka^i\\
=&\nabla([\tau(X_j), A_i]\otimes a^i)+A_i[\tau(X_j),\tau(X_k)]\otimes \partial^ka^i.
\end{align*}
To make the following computation more clear, let us denote:
\begin{align*}
\alpha:=&\nabla([\tau(X_j), A_i]\otimes a^i),\\
\beta:=&A_i[\tau(X_j),\tau(X_k)]\otimes \partial^ka^i.
\end{align*}

First we study $\alpha$. Since $A_i\otimes a^i\in \mathcal{C}_{\tau}(\mathfrak{g})$, by Proposition \ref{criforcfa} we have:
$$
\alpha=\nabla([\tau(X_j), A_i]\otimes a^i)=\nabla(A_i\otimes\{a^i, X_j\})
$$
From the definition of the Poisson bracket on $S(\mathfrak{g})$, we know that
$$
\{a^i, X_j\}=c^r_{sl}X_r\cdot  \partial^sa^i\cdot \partial^lX_j=c^r_{sl}X_r\cdot  \partial^sa^i\cdot \delta^l_j=c^r_{sj}X_r \partial^sa^i.
$$
Therefore
\begin{align*}
\alpha=&\nabla(A_i\otimes c^r_{sj}X_r\cdot  \partial^sa^i)\\
=& A_i\tau(X_l)\otimes \partial^l (c^r_{sj}X_r\cdot  \partial^sa^i)\\
=& A_i\tau(X_l)\otimes c^r_{sj}( \partial^l(X_r)\partial^sa^i+X_r\cdot \partial^l\partial^sa^i)\\
=&A_i\tau(X_l)\otimes c^r_{sj}\delta^l_r\cdot \partial^sa^i+A_i\tau(X_l)\otimes c^r_{sj}X_r\cdot \partial^l\partial^sa^i\\
=&A_i\tau(X_r)\otimes c^r_{sj}\partial^sa^i+A_i\tau(X_l)\otimes c^r_{sj}X_r\cdot \partial^s\partial^la^i.
\end{align*}
Nevertheless, we have
$$
c^r_{sj}X_r\cdot \partial^s\partial^la^i=\{\partial^la^i, X_j\}
$$
As a result
\begin{equation}\label{I}
\alpha=A_i\tau(X_r)\otimes c^r_{sj}\partial^sa^i+A_i\tau(X_l)\otimes \{\partial^la^i, X_j\}
\end{equation}

As for $\beta$, we know
\begin{align*}
\beta=&A_i[\tau(X_j),\tau(X_k)]\otimes \partial^ka^i\\
=&A_i\tau([X_j,X_k])\otimes \partial^ka^i.
\end{align*}
We know that $[X_j,X_k]=c^r_{jk}X_r$ hence
$$
\tau([X_j,X_k])=\tau(c^r_{jk}X_r)=c^r_{jk}\tau(X_r)=-c^r_{kj}\tau(X_r).
$$

As a result
\begin{equation}\label{II}
\begin{split}
\beta=&A_i\tau([X_j,X_k])\otimes \partial^ka^i\\
=&-A_ic^r_{kj}\tau(X_r)\otimes \partial^ka^i\\
=&-A_i\tau(X_r)\otimes c^r_{kj}\partial^ka^i.
\end{split}
\end{equation}

Put Equation (\ref{I}) and (\ref{II}) together we get
\begin{align*}
&\text{the left hand side of }(\ref{eq of image})\\
=&\alpha+\beta\\
=&A_i\tau(X_r)\otimes c^r_{sj}\partial^sa^i+A_i\tau(X_l)\otimes \{\partial^la^i, X_j\}-A_i\tau(X_r)\otimes c^r_{kj}\partial^ka^i\\
=&A_i\tau(X_l)\otimes \{\partial^la^i, X_j\}\\
=&\text{ the right hand side of }(\ref{eq of image}).
\end{align*}
This finishes the proof. \end{proof}

Now with Proposition \ref{nabla maps cfa to cfa}, we can say that the twisted gradient map $\nabla$ is a $\mathbb{C}$-linear map from $\mathcal{C}_{\tau}(\mathfrak{g})$ to $\mathcal{C}_{\tau}(\mathfrak{g})$. In other words, $\nabla$ can be considered as a Hochschild $1$-cochain. see Appendix \ref{app:Hochschild cohomology} for a review of Hochschild cohomology.

\begin{rmk}
In general $\nabla$ is not a Hochschild $1$-cocycle, see Theorem \ref{poisson bracket is a coboundary} below.
\end{rmk}

Before moving on to the next section, we give a direct application of the map $\nabla$.

\begin{coro}[\cite{Ki1} Section 1]\label{famili algebra not only scalar}
When the Lie algebra $\mathfrak{g}$ is semisimple and $\tau$ is a nontrivial irreducible representation, the classical family algebra $\mathcal{C}_{\tau}(\mathfrak{g})$ is more than $I(\mathfrak{g})$, i.e. $I(\mathfrak{g})\varsubsetneq \mathcal{C}_{\tau}(\mathfrak{g})$, and we also have $Z(\mathfrak{g})\varsubsetneq \mathcal{Q}_{\tau}(\mathfrak{g})$.
\end{coro}
\begin{proof} Let Cas be the quadratic Casimir element in $I(\mathfrak{g})$, $\deg \text{Cas}=2$. Then by Proposition \ref{nabla maps cfa to cfa}, we know that $\nabla( \text{Cas})\in \mathcal{C}_{\tau}(\mathfrak{g})$ but $\deg \nabla (\text{Cas})=1$. Since $\tau$ is nontrivial we know that $\nabla (\text{Cas})\neq 0$. On the other hand, since $\mathfrak{g}$ is semisimple, there is no nonzero degree-$1$ element in $I(\mathfrak{g})$, therefore $\nabla (\text{Cas})\notin I(\mathfrak{g})$ hence $I(\mathfrak{g})\varsubsetneq \mathcal{C}_{\tau}(\mathfrak{g})$.

Since there is a PBW map $\FPBW: \mathcal{C}_{\tau}(\mathfrak{g})\rightarrow \mathcal{Q}_{\tau}(\mathfrak{g})$ which maps $I(\mathfrak{g})$ to $Z(\mathfrak{g})$, we know that $Z(\mathfrak{g})\varsubsetneq \mathcal{Q}_{\tau}(\mathfrak{g})$. \end{proof}

\begin{rmk}
In fact, in Example \ref{sl2 standard rep}, the element $M$ is obtained in the same way as $\nabla (\text{Cas})$ in the above corollary.
\end{rmk}

\begin{rmk}
The definition of  $\nabla$ is motivated by the construction of the element $M_P$ defined in Section 1 of \cite{Ki1}. Nevertheless in that paper $M_P$ is defined only for $P\in I(\mathfrak{g})$ and here we extend the domain to all $\mathcal{C}_{\tau}(\mathfrak{g})$.
\end{rmk}

\subsection{The relation between the twisted gradient map and the Poisson bracket}
In this subsection, we build up the relation between $\nabla$ and the Poisson bracket $P$.

First we review some notations of Hochshchild cohomology. Notice that $\nabla: \mathcal{C}_{\tau}(\mathfrak{g})\rightarrow \mathcal{C}_{\tau}(\mathfrak{g})$ is a Hochshchild $1$-cochain, i.e.
$$
\nabla\in C^1(\mathcal{C}_{\tau}(\mathfrak{g}), \mathcal{C}_{\tau}(\mathfrak{g})).
$$

Let
$$
\dH : C^1(\mathcal{C}_{\tau}(\mathfrak{g}), \mathcal{C}_{\tau}(\mathfrak{g}))\rightarrow C^2(\mathcal{C}_{\tau}(\mathfrak{g}), \mathcal{C}_{\tau}(\mathfrak{g}))
$$
be the differential map in the Hochschild complex.

Let $\mathcal{A}, \mathcal{B}\in \mathcal{C}_{\tau}(\mathfrak{g})$. Then by the definition of $\dH$, we have
\begin{equation}
(\dH \nabla)(\mathcal{A},\mathcal{B})=\mathcal{A}\nabla(\mathcal{B})-\nabla(\mathcal{A}\mathcal{B})+\nabla(\mathcal{A})\mathcal{B}.
\end{equation}

The following theorem is the main result of this paper.
\begin{thm}\label{poisson bracket is a coboundary}
For any $\mathcal{A}=A_i\otimes a^i,~\mathcal{B}=B_j\otimes b^j\in \mathcal{C}_{\tau}(\mathfrak{g})$, we have
\begin{equation}\label{poisson to d 1}
\{A,B\}=-\mathcal{A}\nabla(\mathcal{B})+\nabla(\mathcal{A}\mathcal{B})-\nabla(\mathcal{A})\mathcal{B}.
\end{equation}
In other words
\begin{equation}\label{poisson to d 2}
P+\dH \nabla=0
\end{equation}
as elements in the Hochschild $2$-cochain $C^2(\mathcal{C}_{\tau}(\mathfrak{g}), \mathcal{C}_{\tau}(\mathfrak{g}))$. Therefore the Poisson bracket is a Hochschild coboundary in $C^2(\mathcal{C}_{\tau}(\mathfrak{g}), \mathcal{C}_{\tau}(\mathfrak{g}))$.
\end{thm}
\begin{proof}
First let us see what is $\nabla(\mathcal{AB})$:
\begin{equation}
\begin{split}
\nabla(\mathcal{AB})=&\nabla(A_iB_j\otimes a^ib^j)\\
=&A_iB_j\tau(X_k)\otimes \partial^k(a^ib^j)\\
=&A_iB_j\tau(X_k)\otimes (\partial^k a^i)b^j+A_iB_j\tau(X_k)\otimes a^i(\partial^k b^j)
\end{split}
\end{equation}
To make the computation more clear, let us denote:
\begin{align*}
\xi:=&A_iB_j\tau(X_k)\otimes (\partial^k a^i)b^j,\\
\eta:=&A_iB_j\tau(X_k)\otimes a^i(\partial^k b^j).
\end{align*}
Then
\begin{equation}\label{nabla(ab)}
\nabla(AB)=\xi+\eta.
\end{equation}

It is easy to see that $\eta= \mathcal{A}\nabla(\mathcal{B})$. In fact
\begin{equation}\label{anablab}
\eta=A_iB_j\tau(X_k)\otimes a^i(\partial^k b^j)=(A_i\otimes a^i)\cdot (B_j\tau(X_k)\otimes\partial^k b^j)=\mathcal{A}\nabla(\mathcal{B}).
\end{equation}

On the other hand, $\xi\neq (\nabla\mathcal{ A})\mathcal{B}$ in general. We know that
$$
\xi=A_iB_j\tau(X_k)\otimes (\partial^k a^i)b^j
$$
and
\begin{align*}
(\nabla\mathcal{ A})\mathcal{B}=&(A_i\tau(X_k)\otimes \partial^k a^i)\cdot (B_j\otimes b^j)\\
=& A_i\tau(X_k)B_j\otimes (\partial^k a^i)b^j.
\end{align*}
Therefore
\begin{equation}\label{I-nabla}
\begin{split}
\xi-(\nabla\mathcal{ A})\mathcal{B}=&(A_iB_j\tau(X_k)-A_i\tau(X_k)B_j)\otimes (\partial^k a^i)b^j\\
=&A_i[B_j, \tau(X_k)]\otimes (\partial^k a^i)b^j.
\end{split}
\end{equation}

We need to further simplify  the expression $A_i[B_j,\tau(X_k)]\otimes (\partial^k a^i)b^j$. In fact we have the following lemma
\begin{lemma}\label{lemma for deriv to Poisson bracket}
any $\mathcal{A}=A_i\otimes a^i,~\mathcal{B}=B_j\otimes b^j\in \mathcal{C}_{\tau}(\mathfrak{g})$, we have
\begin{equation}
A_i[B_j, \tau(X_k)]\otimes (\partial^k a^i)b^j=A_iB_j\otimes\{a^i, b^j\}=\{\mathcal{A},\mathcal{B}\}.
\end{equation}
\end{lemma}

\begin{proof}[Proof of Lemma \ref{lemma for deriv to Poisson bracket}] First by Lemma \ref{lemma: commute matrices and Lie terms} we have
$$
A_i[B_j,\tau(X_k)]\otimes (\partial^k a^i)b^j=(A_i\otimes \partial^k a^i)\cdot ([B_j,\tau(X_k)]\otimes b^j).
$$

Since $\mathcal{B}=B_j \otimes b^j$ \textbf{is contained in} $\mathcal{C}_{\tau}(\mathfrak{g})$, by Proposition \ref{criforcfa} we know that
\begin{align*}
&(A_i\otimes \partial^k a^i)\cdot (B_j\otimes \{X_k, b^j\})\\
=& A_iB_j\otimes (\partial^k a^i)\cdot \{X_k, b^j\}\\
=& A_iB_j\otimes\{a^i, b^j\}\\
=& \{\mathcal{A},\mathcal{B}\}.
\end{align*}

This proves Lemma \ref{lemma for deriv to Poisson bracket}.
\end{proof}

By Lemma \ref{lemma for deriv to Poisson bracket} and Equation (\ref{I-nabla}) we have
\begin{equation}\label{nablaab}
\xi=\{\mathcal{A},\mathcal{B}\}+(\nabla \mathcal{A})\mathcal{B}.
\end{equation}

Put equations (\ref{nabla(ab)}), (\ref{anablab}) and (\ref{nablaab}) together, we have:
\begin{equation}
\begin{split}\label{dnablatoPoisson}
&\nabla(\mathcal{AB})-\mathcal{A}\nabla(\mathcal{B})-\nabla(\mathcal{A})\mathcal{B}\\
=&\xi+\eta-\mathcal{A}\nabla(\mathcal{B})-\nabla(\mathcal{A})\mathcal{B}\\
=&\{\mathcal{A},\mathcal{B}\}+(\nabla \mathcal{A})\mathcal{B}+\mathcal{A}\nabla(\mathcal{B})-\mathcal{A}\nabla(\mathcal{B})-(\nabla \mathcal{A})\mathcal{B}\\
=&\{\mathcal{A},\mathcal{B}\}.
\end{split}
\end{equation}

This finishes the proof of Theorem \ref{poisson bracket is a coboundary}.
\end{proof}

\begin{ctn}
Although both the twisted gradient map $\nabla$ and the Poisson bracket $P$ can be defined on the larger algebra $\widetilde{\mathcal{C}_{\tau}}(\mathfrak{g})$, we do \textbf{not} have the relation
$$
\{\mathcal{A},\mathcal{B}\}=-\dH \nabla(\mathcal{A},\mathcal{B})
$$
for any $\mathcal{A},\mathcal{B}\in \widetilde{\mathcal{C}_{\tau}}(\mathfrak{g})$. Actually in the proof of Lemma \ref{lemma for deriv to Poisson bracket} we see that it is necessary to have $\mathcal{B}\in \mathcal{C}_{\tau}(\mathfrak{g})$.
\end{ctn}

From the view point of deformation theory (Proposition \ref{infinites trivial deform}), we have the following corollary.

\begin{coro}\label{Poisson family algebra infini trivial deform}
The  $1$-parameter formal deformation from $\mathcal{C}_{\tau}(\mathfrak{g})$ to $\mathcal{Q}_{\tau}(\mathfrak{g})$ is infinitesimally trivial.
\end{coro}
\begin{proof}
 We know in Proposition \ref{first terms of the star product on family algebras} that the first order  $1$-parameter formal deformation $m_1$ is $\frac{1}{2}P$, therefore this corollary is just a direct consequence of Theorem \ref{poisson bracket is a coboundary}.
\end{proof}

Using Theorem \ref{poisson bracket is a coboundary} we can also give an alternative proof of Proposition \ref{poisson bracket mas cfa to cfa} as follows.

\begin{coro}[Proposition \ref{poisson bracket mas cfa to cfa}]\label{coro: poisson bracket mas cfa to cfa}
For any $\mathcal{A}, \mathcal{B}\in \mathcal{C}_{\tau}(\mathfrak{g})$, we have that $\{\mathcal{A}, \mathcal{B}\}$ is still in $\mathcal{C}_{\tau}(\mathfrak{g})$.
\end{coro}
\begin{proof}
 In the proof of Theorem \ref{poisson bracket is a coboundary}, we do not require a priori that  $\{\mathcal{A}, \mathcal{B}\}\in \mathcal{C}_{\tau}(\mathfrak{g})$. Now by Proposition \ref{nabla maps cfa to cfa} we know that $\nabla$ maps $\mathcal{C}_{\tau}(\mathfrak{g})$ to $\mathcal{C}_{\tau}(\mathfrak{g})$ and
 from Theorem \ref{poisson bracket is a coboundary} we also know that
 $$
\{\mathcal{A},\mathcal{B}\}=-\dH \nabla(\mathcal{A},\mathcal{B})
$$
hence we get the result we want.
 \end{proof}

\subsection{Digression: An alternative of the twisted gradient map}
 In this subsection we want to show that the twisted gradient map $\nabla$ defined in Definition \ref{defofnabla} is NOT the unique map which satisfies
  $$
\{\mathcal{A},\mathcal{B}\}=-\dH \nabla(\mathcal{A},\mathcal{B}).
$$

In fact we define a map $\nabla':\widetilde{\mathcal{C}_{\tau}}(\mathfrak{g})\rightarrow \widetilde{\mathcal{C}_{\tau}}(\mathfrak{g})$ to be
\begin{equation}
\nabla'(A_i\otimes a^i):=\tau(X_k)A_i\otimes \partial^ka^i.
\end{equation}
Similar to Proposition \ref{nabla maps cfa to cfa}, we can check that $\nabla'$ also maps $\mathcal{C}_{\tau}(\mathfrak{g})$ to $\mathcal{C}_{\tau}(\mathfrak{g})$.

\begin{rmk}
The difference between the definition of $\nabla$ and $\nabla'$ is: for $\nabla$, the matrix $\tau(X_k)$ is multiplied from the right; while  for $\nabla'$, the matrix $\tau(X_k)$ is multiplied from the left.
\end{rmk}

In general $\nabla' \neq \nabla $ and we want to find their difference. First we define the \emph{first Chern class} on $\widetilde{\mathcal{C}_{\tau}}(\mathfrak{g})$ following \cite{CalRoss2011} Section 1.1.

\begin{defi}\label{first chern class in g}
The first Chern class $c_1$ is a map $\widetilde{\mathcal{C}_{\tau}}(\mathfrak{g})\rightarrow \widetilde{\mathcal{C}_{\tau}}(\mathfrak{g})$, $c_1:=\tr(\ad)$. More precisely
\begin{equation}
\begin{split}
c_1: \widetilde{\mathcal{C}_{\tau}}(\mathfrak{g})&\longrightarrow \widetilde{\mathcal{C}_{\tau}}(\mathfrak{g})\\
A\otimes a &\mapsto A\otimes c^j_{ij}\partial^i a.
\end{split}
\end{equation}

It is easy to check that $c_1$ is $\mathfrak{g}$-invariant hence $c_1$ maps $\mathcal{C}_{\tau}(\mathfrak{g})$ to $\mathcal{C}_{\tau}(\mathfrak{g})$. Moreover, it is also easy to check that the first Chern class is closed in the Hochschild cochain. In other words, $c_1\in Z^1(\widetilde{\mathcal{C}_{\tau}}(\mathfrak{g}),\widetilde{\mathcal{C}_{\tau}}(\mathfrak{g}))$ and $c_1|_{\mathcal{C}_{\tau}(\mathfrak{g})}\in Z^1(\mathcal{C}_{\tau}(\mathfrak{g}),\mathcal{C}_{\tau}(\mathfrak{g}))$. For simplicity we also write $c_1$ for the restriction $c_1|_{\mathcal{C}_{\tau}(\mathfrak{g})}$.
\end{defi}

Having the first Chern class, we can express the difference between $\nabla$ and $\nabla'$ in $\mathcal{C}_{\tau}(\mathfrak{g})$:

\begin{prop}
In the classical family algebra $\mathcal{C}_{\tau}(\mathfrak{g})$ we have
\begin{equation}\nabla-\nabla'=-c_1.
\end{equation}
\end{prop}
\begin{proof}
 For any $ A_i\otimes a^i\in \mathcal{C}_{\tau}(\mathfrak{g})$
\begin{align*}
\nabla(A_i\otimes a^i)-\nabla'(A_i\otimes a^i)=&A_i\tau(X_k)\otimes \partial^ka^i-\tau(X_k)A_i\otimes \partial^ka^i\\
=& [A_i, \tau(X_k)]\otimes \partial^ka^i\\
=&\partial^k([A_i, \tau(X_k)]\otimes a^i)~(\text{We can move the partial derivative out}).
\end{align*}

Since $ A_i\otimes a^i\in \mathcal{C}_{\tau}(\mathfrak{g})$, we have
$$
[A_i, \tau(X_k)]\otimes a^i= A_i\otimes \{X_k, a^i\}.
$$

Therefore
\begin{align*}
\nabla(A_i\otimes a^i)-\nabla'(A_i\otimes a^i)=&\partial^k(A_i\otimes \{X_k, a^i\})\\
=& A_i\otimes \partial^k( \{X_k, a^i\})\\
=& A_i\otimes \partial^k(c_{kj}^lX_l\cdot \partial^ja^i)\\
=& A_i\otimes (c_{kj}^k\partial^ja^i+c_{kj}^lX_l\cdot \partial^k\partial^ja^i)\\
=& A_i\otimes c_{kj}^k\partial^ja^i+A_i\otimes c_{kj}^lX_l\cdot \partial^k\partial^ja^i\\
=&-c_1(A_i\otimes a^i)+A_i\otimes c_{kj}^lX_l\cdot \partial^k\partial^ja^i
\end{align*}
Since ~$c_{kj}^l$~is anti symmetric with respect to ~$k,j$,~it is easy to see that
$$
A_i\otimes c_{kj}^lX_l\cdot \partial^k\partial^ja^i=0
$$
Hence we get
$$
\nabla(A_i\otimes a^i)-\nabla'(A_i\otimes a^i)=-c_1(A_i\otimes a^i).
$$
\end{proof}

The next corollary tells us that we can replace $\nabla$  by $\nabla'$ in Theorem \ref{poisson bracket is a coboundary}.

\begin{coro}\label{nabla' poisson bracket}
In $\mathcal{C}_{\tau}(\mathfrak{g})$ we have $$\dH \nabla'=\dH \nabla=P, \text{ the Poisson bracket on }\mathcal{C}_{\tau}(\mathfrak{g}).$$
Therefore we can replace $\nabla$  by $\nabla'$ in Theorem \ref{poisson bracket is a coboundary}.
\end{coro}
\begin{proof}
We know that $\dH \nabla=P$ and $\nabla-\nabla'=-c_1$. In Definition \ref{first chern class in g} we also know that $c_1$ is closed, i.e. $\dH\, c_1=0$. Hence we get this corollary.
\end{proof}

At the end of this subsection we should point out that although $\nabla' \neq \nabla $ in general, they are equal in some important cases. Actually  we have the following result.

\begin{prop}\label{semisimple nabla'}
When $\mathfrak{g}$ is a semisimple Lie algebra, we have $\nabla'=\nabla$ in $\mathcal{C}_{\tau}(\mathfrak{g})$.
\end{prop}
\begin{proof}  We know that for semisimle Lie algebra, the adjoint representation is traceless, in other words
$$
c^i_{ij}=0 \text{ for any } j.
$$
Therefore $c_1=0$ for semisimple $\mathfrak{g}$, hence we get the result. \end{proof}


\section{Further topics: Mackey's analogue and The quantization of the family algebras}\label{section: quantization of family algebras}
This section is a survey of further topics. In this section we restrict to the case that $\mathfrak{g}$ is a complex semisimple Lie algebra and the representation $\tau$ to be a simple representation of $\mathfrak{g}$.\\

In 1975 G. Mackey (\cite{Mackey1975}) studied the analogies between the representations of a semisimple Lie group  $G$  and those of its Cartan motion group $G_c$. Later N. Higson find the relation between family algebras and Mackey's analogue in \cite{Higson2011}  as we have mentioned in the introduction.

\begin{rmk}
In fact, Higson introduced the \emph{spherical Hecke algebras} $\mathcal{R}(\mathfrak{g},\tau)$ and  $\mathcal{R}(\mathfrak{g}_c,\tau)$ respectively. These algebras have the importance that the irreducible $\mathcal{R}(\mathfrak{g},\tau)$ modules are $1$-$1$ correspondent to irreducible $(\mathfrak{g},K)$-modules of $G$ with nonzero $\tau$-isotypical component, and the similar result holds for $\mathcal{R}(\mathfrak{g}_c,\tau)$. When $G$ is complex semisimple, Higson proved that the spherical Hecke algebras are isomorphic to the family algebras. For details see \cite{Higson2011}.
\end{rmk}

Let $\mathfrak{h}\subset \mathfrak{g}$ be the Cartan subalgebra. Higson also constructed the \emph{generalized Harish-Chandra homomorphisms}:
\begin{equation}
\begin{split}
&\text{GHC}_{\tau}: \mathcal{R}(\mathfrak{g},\tau)\rightarrow U(\mathfrak{h}) \\
&\text{GHC}_{\tau,c}: \mathcal{R}(\mathfrak{g}_c,\tau)\rightarrow S(\mathfrak{h})
\end{split}
\end{equation}
and relates them to the admissible duals of $G$ and $G_c$ with minimal $K$-type $\tau$.

The Mackey's analogue for admissible dual of complex semisimple $G$ has the following form:
\begin{thm}[\cite{Higson2011}, Section 8]
Under the identification $U(\mathfrak{h})\cong S(\mathfrak{h})$, the two homomorphisms $\text{GHC}_{\tau}$ and $\text{GHC}_{\tau,c}$ has the same image.
\end{thm}

In the end of \cite{Higson2011}, Higson proposed the problem of constructing a \emph{quantization map} $Q$ between $\mathcal{C}_{\tau}(\mathfrak{g})$ and $\mathcal{Q}_{\tau}(\mathfrak{g})$ such that the following diagram commutes.
\begin{equation}
\begin{tikzpicture}[node distance=3cm, auto]
  \node (A1) {$\mathcal{C}_{\tau}(\mathfrak{g})$};
  \node [right of=A1] (A2) {$\mathcal{Q}_{\tau}(\mathfrak{g})$};
  \node [below of=A1, node distance=1.5cm] (B1) {$S(\mathfrak{h})$};
  \node [below of=A2, node distance=1.5cm] (B2) {$U(\mathfrak{h})$};

\draw[->] (A1) to node {$\text{\tiny{GHC}}_{\tau,c}$} (B1);
\draw[->,dashed] (A1) to node {$Q$}(A2);
\draw[->] (A2) to node {$\text{\tiny {GHC}}_{\tau}$}(B2);
\draw[->] (B1) to node {$\cong$}(B2);
\end{tikzpicture}
\end{equation}
Here $Q$ is a vector space isomorphism but need not to be an algebraic isomorphism.

\begin{rmk}
According to Theorem \ref{poisson bracket is a coboundary} and Corollary \ref{Poisson family algebra infini trivial deform}, the  $1$-parameter formal deformation from $\mathcal{C}_{\tau}(\mathfrak{g})$ to $\mathcal{Q}_{\tau}(\mathfrak{g})$ is infinitesimally trivial, which suggests that there exists strong relations between them. Actually Theorem  \ref{poisson bracket is a coboundary} can be considered as the first step in the solution of the quantization problem of the family algebras.
\end{rmk}

\section*{Acknowledgement} Z.W. wishes to thank Alexandre Kirillov, Vasily Dolgushev, Eckhard Meinrenken, Valery Lunts and in particular Nigel Higson for helpful discussions and comments.

This paper is part of the doctoral dissertation of Z.W. at University of Pennsylvania and Z.W would like to thank his doctoral advisor Jonathan Block for his help on this topic and encouragements.

\appendix
\numberwithin{equation}{section}
\renewcommand\thethm{\Alph{section}.\arabic{thm}}
\renewcommand\thedefi{\Alph{section}.\arabic{defi}}
\renewcommand\thermk{\Alph{section}.\arabic{rmk}}
\renewcommand\theeg{\Alph{section}.\arabic{eg}}
\renewcommand\thectn{\Alph{section}.\arabic{ctn}}
\section{Hochschild cohomology}\label{app:Hochschild cohomology}
Let us review the theory of Hochschild cohomology in this appendix. For reference see \cite{Weibel1994} or \cite{CalRoss2011} Section 2.

Let $A$ be an associative $\mathbb{C}$-algebra. The  \emph{Hochschild complex} $C^{\bullet}(A,A)$ is defined as follows:
\begin{equation}
C^{n}(A,A):= \Hom_{\mathbb{C}}(A^{\otimes n},A),~n\geqslant 0.
\end{equation}
The differential $\dH$ is defined on homogeneous elements $f\in C^{n}(A,A)$ by the formula
\begin{equation}\label{defi of dH}
\begin{split}
(\dH(f))(a_0,a_1,\ldots,a_n):= &a_0 f(a_1,\ldots,a_n)+\sum_{k=1}^n(-1)^kf(a_0,\ldots, a_{k-1}a_k,\ldots,a_n)\\
+&(-1)^{n+1}f(a_0,\ldots,a_{n-1})a_n.
\end{split}
\end{equation}
We see that $\dH f\in C^{n+1}(A,A)$. Actually we can prove that $\dH\circ \dH=0$ therefore  $C^{\bullet}(A,A)$ is a cochain complex.

The \emph{Hochschild cohomology} of $A$ is defined as the cohomology group of the cochain complex $C^{\bullet}(A,A)$, and we denote it by $\HH^{\bullet}(A,A)$ or for short $\HH^{\bullet}(A)$:
\begin{equation}
\HH^{n}(A):= \text{H}^n(C^{\bullet}(A,A)).
\end{equation}

Now let us look at the case $n=2$. The following result is easy to get:
\begin{prop}\label{hochschild 2 coboundary and cocycle}
Let $f\in C^{2}(A,A)=\Hom_{\mathbb{C}}(A\otimes A,A)$. Then $f$ is a $2$-coboundary if and only if there exists a $g\in C^{1}(A,A)=\Hom_{\mathbb{C}}(A,A)$ such that for any $a,b\in A$
\begin{equation}
f(a,b)=ag(b)-g(ab)+g(a)b.
\end{equation}

Moreover, $f$ is a $2$-cocycle  if and only if for any $a,b, \in A$
\begin{equation}
af(b,c)-f(ab,c)+f(a,bc)-f(a,b)c=0.
\end{equation}
\end{prop}
\begin{proof} Direct check by definition. \end{proof}

\section{The Gerstenhaber bracket on Hochschild cochains and cohomologies}\label{app:HH and Gerstenhaber}
In this section we give a quick review of the \emph{Gerstenhaber bracket}. For more details and proofs see \cite{Gerstenhaber1963} or \cite{BlockGetzler1992} Section 1. For further topics see the survey \cite{DTT2009}.

First, we define an operation $\circ: C^k(A,A)\otimes C^l(A,A)\rightarrow C^{k+l-1}(A,A)$. Let $f_1\in C^k(A,A)$ and $f_2\in C^l(A,A)$,
\begin{equation}\label{define of circ on HC}
\begin{split}
(f_1\circ f_2)&(a_1,\ldots,a_{k+l-1}):=\\
=&\sum_{i=0}^{k-1}(-1)^{(k-i-1)(l-1)}f_1(a_1\ldots,a_i,f_2(a_{i+1},\ldots, a_{i+l}),a_{i+l+1},\ldots,a_{k+l-1}).
\end{split}
\end{equation}

In particular, for $2$-cochains we have
\begin{prop}\label{circ on 2-cochain}
Let $f_1, f_2\in C^2(A,A)$, then $f_1\circ f_2\in C^3(A,A)$ and is given by
\begin{equation}\label{define of circ on 2 cochain}
(f_1\circ f_2)(a_1,a_2,a_3)=f_1(f_2(a_1,a_2),a_3)-f_1(a_1,f_2(a_2,a_3)).
\end{equation}
In particular, for $f\in C^2(A,A)$ we have
\begin{equation}\label{define of circ with itself on 2 cochain}
(f \circ f )(a_1,a_2,a_3)=f(f(a_1,a_2),a_3)-f(a_1,f(a_2,a_3)).
\end{equation}
\end{prop}
\begin{proof} This is exactly the definition. \end{proof}

The \emph{Gerstenhaber bracket} is defined to be
\begin{equation}\label{define of Gerstenhaber bracket on HC}
[f_1, f_2]_{\text{G}}:= f_1\circ f_2-(-1)^{(k-1)(l-1)}f_2\circ f_1.
\end{equation}

The Gerstenhaber bracket is a Lie bracket. In fact we have the following

\begin{thm}\label{gerstenhaber bracket is a Lie bracket}
The operation "$\,\circ$" gives a pre-Lie algebra structure on  $C^{\bullet-1}(A,A)$. Therefore we obtain that
$(C^{\bullet-1}(A,A), [\,,\,]_{\text{G}})$ is a graded Lie algebra.
\end{thm}
\begin{proof} See \cite{Gerstenhaber1963}. \end{proof}

\begin{prop}\label{gerstenhaber with itself}
Let $f\in C^2(A,A)$, then
\begin{equation}\label{gerstenhaber with itself on 2-cochain}
[f, f]_{\text{G}}= 2f\circ f.
\end{equation}
\end{prop}
\begin{proof}  We get this directly from the definitions.  \end{proof}

In fact $\dH$ is an inner derivation under the Gerstenhaber bracket. More precisely, let
$\mu: A\otimes A\rightarrow A$ be the multiplication map in $A$. Then $\mu\in C^2(A,A)$ and we have the following
\begin{prop}\label{dH is comm of mu}
For any $f\in C^k(A,A)$, we have
\begin{equation}
\dH f= [\mu, f]_{\text{G}}\in C^{k+1}(A,A).
\end{equation}
We also have $[\mu, \mu]_{\text{G}}=0$.
\end{prop}
\begin{proof}  Compare the definition of $\dH$ in Equation (\ref{defi of dH}) and the definition of the Gerstenhaber bracket in  Equation (\ref{define of circ on HC}) and Equation (\ref{define of Gerstenhaber bracket on HC}). The fact that $[\mu, \mu]_{\text{G}}=0$ is exactly the associativity of $\mu$. \end{proof}

As a result, we have the following theorem:
\begin{thm}\label{gerstenhaber bracket reduce to HH}
The Gerstenhaber bracket is compatible with the Hochschild differential $\dH$. In other words, for any $f_1\in C^k(A,A)$ and $f_2\in C^l(A,A)$, we have
\begin{equation}\label{dH and Gerstenhaber Leibniz}
\dH ([f_1, f_2]_{\text{G}}) =[\dH f_1, f_2]_{\text{G}}+(-1)^{k-1}[f_1,\dH f_2]_{\text{G}}.
\end{equation}

Therefore the Gerstenhaber bracket reduces to the Hochschild cohomology $\HH^{\,\bullet-1}(A)$.
\end{thm}
\begin{proof} Since $\dH$ is an inner derivation according to Propostion \ref{dH is comm of mu}, Equation (\ref{dH and Gerstenhaber Leibniz}) is a consequence of  the graded-Jacobi identity of the graded Lie algebra $(C^{\bullet-1}(A,A), [\,,\,]_{\text{G}})$.\end{proof}

\section{$\HH^{\bullet}(A)$ and the deformations of $A$}\label{app:HH and deformation}
The Hochschild cohomology plays an important role in the deformation theory. Let us summarize some results in the deformation theory of algebras in this appendix. For more details see \cite{Gerstenhaber1964}.

Let $A$ be an associative $\mathbb{C}$ algebra (in fact we can replace $\mathbb{C}$ by any field). A deformation of the algebra structure of $A$ means that we fix $A$ as a $\mathbb{C}$-vector space and change the multiplication operation on $A$. Actually there are many kinds of deformations, like analytic, algebraic, formal, global, etc., and in this paper we focus on formal deformation, in particular \emph{$1$-parameter formal deformation} of algebras. In more details let $\mathbb{C}[[t]]$ be the formal power series of $t$ and we define
\begin{equation}
 A[[t]]:=A\otimes_{\mathbb{C}}\mathbb{C}[[t]].
\end{equation}

$A[[t]]$ is obviously a $\mathbb{C}[[t]]$-module.

A $1$-parameter formal deformation of the algebra structure on $A$ is given by a map
\begin{equation}\label{deformation of algebra1}
m : A[[t]]\otimes A[[t]]\longrightarrow A[[t]]
\end{equation}
where $m$ is required to be $\mathbb{C}[[t]]$-bilinear. So we only need to know the value of $m$ on $A\otimes A$. Moreover we require that $
m(a,b)\equiv ab \text{ mod } t$ for  $a,b\in A$.

For any $a,b\in A$, we can write $m(a,b)$ as
\begin{equation}\label{deformation of algebra2}
m(a,b)=ab+\sum_{k=0}^{\infty}t^km_k(a,b).
\end{equation}
We see that each $m_k$ belongs to $C^2(A,A)$.

\begin{rmk} The element $t$ is called the deformation parameter. To get an informal idea of deformation theory, we can evaluate at $t=0$ and get the original multiplication on $A$. On the other hand if we evaluate at $t\neq 0$, omit the convergence problem, we get a new binary operation $A\otimes A\rightarrow A$.
\end{rmk}

\begin{rmk} We can also talk about more general formal deformation of $A$, where the algebra $\mathbb{C}[[t]]$ is replaced by  a \emph{complete local augmented $\mathbb{C}$-algebra}, see \cite{doubek2007deformation} Section 3.

Moreover we also have formal deformation theory of Lie algebras,see \cite{fialowski1988example} for a detailed introduction.
\end{rmk}

 As a multiplication, the map $m$ needs to satisfy the \emph{associativity law}.

\begin{thm}[$1$-parameter formal deformation, see \cite{Gerstenhaber1964} Chapter I.1]\label{when deformation m is associative}
Let $m(a,b)=ab+\sum_{k=0}^{\infty}t^km_k(a,b)$ as in Equation (\ref{deformation of algebra2}). Then $m$ satisfies the associativity law if and only if
for each $k \geqslant 1$, we have
\begin{equation}\label{differential and associativity law}
\dH m_k+\frac{1}{2}\sum_{i=1}^{k-1}[m_i,m_{k-i}]=0.
\end{equation}
If this holds, we say that $m$ gives a formal deformation of $A$.
\end{thm}
\begin{proof}
The associativity law tells us that for any $a,b,c\in A$, we have
\begin{equation}\label{assoc law}
m(a,m(b,c))-m(m(a,b),c)=0.
\end{equation}

Now consider $m$ as an element in $C^2(A[[t]],A[[t]])$, then Equation (\ref{assoc law}) is exactly
\begin{equation}\label{assoc law in Gerstenhaber bracket}
[m,m]_{\text{G}}=0.
\end{equation}

We write $m=\mu+\sum_{k=1}^{\infty}t^km_k$ where $\mu$ is the original multiplication on $A$. Then because we know
$$[\mu,f]_{\text{G}}=\dH f \text{ and } [\mu,\mu]_{\text{G}}=0$$ in Proposition \ref{dH is comm of mu}, Equation (\ref{assoc law in Gerstenhaber bracket}) becomes the \emph{Maurer-Cartan Equation}
\begin{equation}\label{Maurer-Cartan Equation}
\dH (\sum_{k=1}^{\infty}t^km_k)+\frac{1}{2}[\sum_{k=1}^{\infty}t^km_k,\sum_{k=1}^{\infty}t^km_k]_{\text{G}}=0.
\end{equation}

In the expansion of Equation (\ref{Maurer-Cartan Equation}), we take the $t^k$ term and get Equation (\ref{differential and associativity law}).
\end{proof}

\begin{coro}[Infinitesimal deformation]\label{deformation up to 2 and 3}
The map $m$ satisfies the associativity law$\mod t^2$ if and only if $\dH\, m_1=0$, i.e. for any $a,b,c\in A$, we have
\begin{equation}\label{dH m 1=0}
am_1(b,c)-m_1{ab,c}+m_1(a,bc)-m_1(a,b)c=0.
\end{equation}
If this holds, we say that $m$ gives an infinitesimal deformation of $A$.

Moreover, $m$  satisfies the associativity law$\mod t^3$ if and only if $\dH\, m_1=0$ together with
\begin{equation}\label{dH m 2 and Gerstenhaber}
\dH\, m_2+\frac{1}{2}[m_1,m_1]_{\text{G}}=0.
\end{equation}
The above equation is equivalent to
\begin{equation}\label{dH m 2 and circ}
\dH\, m_2+m_1\circ m_1 =0
\end{equation}
\end{coro}
\begin{proof} This is an direct corollary of Theorem \ref{when deformation m is associative}. \end{proof}

On the other hand, we need to know when the $1$-parameter formal deformation $m$ is trivial. In other words, whether or not we can find an algebraic isomorphism
\begin{equation}\label{phi alge isom}
\theta:~(A[[t\,]], \mu)\longrightarrow (A[[t\,]], m)
\end{equation}
where $\theta$ is $\mathbb{C}[[t]]$-linear and is given by
\begin{equation}\label{theta formula}
\theta(a)=a+\sum_{k=1}^{\infty}t^k\theta_k(a).
\end{equation}

The requirement for $\theta$ is for any $a,b\in A$
\begin{equation}\label{the requirement for theta}
\theta(ab)=m(\theta(a),\theta(b)).
\end{equation}

The existence of $\theta$ is a complicated issue. Nevertheless as a first step we have:
\begin{prop}[Infinitesimally trivial deformation]\label{infinites trivial deform}
There exists a $\theta_1\in C^1(A,A)$ such that $\theta=id+t\theta_1$ satisfyies Equation (\ref{the requirement for theta})$\mod t^2$ if and only if $m_1\in B^2(A,A)$. If this holds, we say that $m$ is an infinitesimally trivial deformation of $A$.
\end{prop}
\begin{proof} We expand both sides of Equation (\ref{the requirement for theta}) and look at the $t$ term we get
\begin{equation}
\theta_1(ab)=\theta_1(a)b+a\theta_1(b)+m_1(a,b)
\end{equation}
In other words
\begin{equation}
m_1+\dH\theta_1=0.
\end{equation}
\end{proof}

Further discussion of the triviality of formal deformations involves the concept of \emph{gauge equivalence} of Maurer-Cartan elements, see \cite{Kontsevich2003} Section 1 or \cite{LPV2013} Chapter 13.

\bibliography{familyalgebrabib}{}
\bibliographystyle{plain}

\end{document}